\documentclass[12pt,reqno]{amsart}
\usepackage{amssymb,delarray}
\usepackage{amsfonts}
\usepackage{epsfig}
\usepackage[all]{xy}

\makeindex{}

\newtheorem{lem}
{Lemma}
\newtheorem{prop}
{Proposition}

{Corollary}

{\catcode`\@=11
\gdef\n@te#1#2{\leavevmode\vadjust{%
 {\setbox\z@\hbox to\z@{\strut#1}%
  \setbox\z@\hbox{\raise\dp\strutbox\box\z@}\ht\z@=\z@\dp\z@=\z@%
  #2\box\z@}}}
\gdef\leftnote#1{\n@te{\hss#1\quad}{}}
\gdef\rightnote#1{\n@te{\quad\kern-\leftskip#1\hss}{\moveright\hsize}}
\gdef\?{\FN@\qumark}
\gdef\qumark{\ifx\next"\DN@"##1"{\leftnote{\rm##1}}\else
 \DN@{\leftnote{\rm??}}\fi{\rm??}\next@}}

\begin{document}
\baselineskip=14.pt plus 2pt 

\title[]{A remark on the non-rationality Problem
for generic cubic fourfolds}
\author[Vik.S.~Kulikov]{Vik. S.~Kulikov}
\address{Steklov Mathematical Institute\\
Gubkina str., 8\\
119991 Moscow \\
Russia} \email{kulikov@mi.ras.ru}

\dedicatory{} \subjclass{}
\thanks{The work  was partially supported
by the RFBR  (05-01-00455),   NWO-RFBR 047.011.2004.026
(RFBR 05-02-89000-$NWO_a$), INTAS (05-1000008-7805), and by
RUM1-2692-MO-05. }
\keywords{}
\begin{abstract}
It is proved that the non-rationality of a generic cubic fourfold
follows from a conjecture on the non-decomposability in the direct
sum of non-trivial polarized Hodge structures of the polarized
Hodge structure on transcendental cycles on a projective surface.
\end{abstract}

\maketitle
\setcounter{tocdepth}{2}


\def\st{{\sf st}}

\setcounter{section}{1}

Before to formulate the main statement of the present note, recall
some definitions and statements related to Hodge structures.

Let $\mathcal H=\{ M, H^{p,n-p},Q\}$  be a {\it polarized Hodge
structure of weight} $n$, that is, $M$ is a free $\mathbb Z$-module,
$Q:M\times M\to \mathbb Z$ is a non-degenerate bilinear form, the
vector space $M\otimes \mathbb C=\oplus_{p=a}^b H^{p,n-p}$ is the
direct sum of its complex subspaces $H^{p,n-p}$ for some $a,b\in
\mathbb Z$, $a\leq b$, such that $H^{n-p,p}=\overline {H^{p,n-p}}$
and
$$
\begin{array}{rrll} Q(u,v) & = & (-1)^nQ(v,u); & \\
Q(u,v) & = & 0 & \text{for}\,\, u\in H^{p,n-p},\,\, v\in
H^{p^\prime,n-p^\prime}, \,\, p\neq n-p^\prime.
\end{array}
$$

A polarized  Hodge structure $\mathcal H=\{ M, H^{p,q},Q\}$ is
called {\it unimodular} if $Q$ is an unimodular bilinear form,
that is, if $(e_i)$ is a free basis of $M$, then the determinant
of the matrix $A=(Q(e_i,e_j))$ is equal to $\pm 1$. Let $\mathcal
H_{1}=\{ M_1, H_1^{p,q},Q_1\}$ and $\mathcal H_{2}=\{ M_2,
H_2^{p,q},Q_2\}$ be two polarized Hodge structures of weight $n$.
A $\mathbb Z$-homomorphism $f:M_1\to M_2$ of $\mathbb Z$-modules
is a {\it morphism of  polarized  Hodge structures} if
$Q_2(f(u),f(v))=Q_1(u,v)$ for all $u,v\in M_1$ and $f$ induces a
morphism of Hodge structures, that is, $f_{\mathbb
C}(H_1^{p,q})\subset H_2^{p,q}$ for all $p,q$, where $f_{\mathbb
C}=f\otimes \text{Id}$. Note that if $f:M_1\to M_2$ is a morphism
of polarized Hodge structures, then $f$ is an embedding, since
$Q_1$ and $Q_2$ are non-degenerate bilinear forms. We say that
$\mathcal H_{1}$ is a {\it polarized Hodge substructure} of
$\mathcal H_{2}$ (and will write $\mathcal H_{1}\subset \mathcal
H_{2}$) if the embedding $M_1\subset M_2$ is a morphism of
polarized Hodge structures. A polarized Hodge structure $\mathcal
H=\{ M_1\oplus M_2, H^{p,q},Q_1\oplus Q_2\}$ is called the {\it
direct sum of polarized Hodge structures} $\mathcal H_{1}$ and
$\mathcal H_{2}$ if the canonical embeddings $M_i\subset M_1\oplus
M_2$ are morphisms of polarized Hodge structures for $i=1,2$. In
the case of the direct sum we have $H^{p,q}=H_{1}^{p,q}\oplus
H_{2}^{p,q}$ for all $p,q$. We say that a polarized Hodge
structure $\mathcal H$ is {\it non-decomposable} if it is not
isomorphic to the direct sum of two non-trivial polarized Hodge
structures.

Let $X$ be a smooth projective manifold defined over the field
$\mathbb C$, $\dim_{\mathbb C} X=n$. It is well known that one can
associate to $X$ a polarized Hodge structure $\mathcal H_X=\{ M_X,
H^{p,q},Q\}$ of weight $n$, where $$M_X=H^n(X,\mathbb Z)/\{
torsion\},$$ $H^{p,q}=H^{p,q}_X\subset H^n(X,\mathbb C)$ are the
spaces of harmonic $(p,q)$-forms on $X$, and $Q=Q_X$ is the
restriction to $M_X\subset M_X\otimes \mathbb C\simeq
H^n(X,\mathbb C)$ of the bilinear form
$$Q(\phi,\psi)=\int_X\phi\wedge\psi .$$
It is well known that the lattice $(M_X, Q_X)$ is unimodular. Denote
by $h^{p,q}=h^{p,q}(X)=\dim H^{p,q}_X=\dim H^q(X,\Omega_X^p)$ the
Hodge numbers of $X$, where $\Omega_X^p$ is the sheaf of holomorphic
$p$-forms on $X$.

\begin{lem} Let $f:X\to V$ be a bi-rational morphism of smooth
projective manifolds, $\dim_{\mathbb C}X=\dim_{\mathbb C}V=n$.
Then the polarized Hodge structure $\mathcal H_{X}$ is decomposed
in the direct sum of the polarized Hodge structure $f^*(\mathcal
H_{V})\simeq \mathcal H_{V}$ and a polarized Hodge structure
$\mathcal H_V^{\bot}$.
\end{lem}
\proof  This statement is well-known in the particular case if
$f=\sigma :X\to V$ is the monoidal transformation with non-singular
center $C\subset V$ (\cite{A}, \cite{D}, \cite{K}; see also
\cite{C-G}). By induction, the correctness of lemma follows for a
composition $f=\sigma_m\circ\dots\circ\sigma_1 :X\to V$ of monoidal
transformations $\sigma_i: X_i\to X_{i-1}$ with non-singular centers
$\widetilde C_{i-1}\subset X_{i-1}$, where $X_0=V$ and $X_m=Y$.

Now, let  $f$ be an arbitrary bi-rational morphism.  Since $\deg
f=1$, then
$$\int_V\phi\wedge\psi =\int_Xf^*(\phi)\wedge f^*(\psi), $$ that is,
$f^*$ is an embedding of the lattice $M_V=H^4(V,\mathbb Z)$ in the
lattice $M_X=H^4(V,\mathbb Z)$. Therefore $f^*: \mathcal H_V\to
\mathcal H_X$ is an embedding of the polarized Hodge structure.
Since the lattices $f^*(M_V)$ and $M_X$ are unimodular, there is a
sublattice $M_V^\bot$ of $M_X$ such that $M_X=f^*(M_V)\oplus
M_V^\bot$ is the direct sum of lattices.

Let us show that this decomposition in the direct sum of lattices
induces a decomposition $\mathcal H_X=f^*(\mathcal H_V)\oplus
\mathcal H_V^\bot$ in the direct sum of polarized Hodge
structures.  For this, since the image of Hodge structure under a
morphism of Hodge structures is a Hodge structure, it suffices to
show that the natural projection $\text{pr}^\bot:M_X\otimes
\mathbb C\to M_V^\bot\otimes \mathbb C$ is a morphism of Hodge
structures. To show that $\text{pr}^\bot$ is a morphism of Hodge
structures, let us consider the bi-rational map $f^{-1}:
V\dashrightarrow X$. By Hironaka Theorem, there is a composition
$\sigma =\sigma_m\circ\dots\circ\sigma_1 :Y\to V$ of monoidal
transformations $\sigma_i: Y_i\to Y_{i-1}$ with non-singular
centers $\widetilde C_{i-1}\subset Y_{i-1}$, where $Y_0=V$ and
$Y_m=Y$, such that $g=f^{-1}\circ \sigma :Y\to X$ is a bi-rational
morphism. We obtain a decomposition
$$\mathcal H_Y=\sigma^*(\mathcal H_{V})\oplus\bigoplus_{i=0}^{m-1}\widetilde{\mathcal
H}_i$$ in the direct sum of polarized Hodge structures, where
$\widetilde{\mathcal H}_i$ is the contribution in $\mathcal H_Y$
of the monoidal transformation $\sigma_{i+1}$. Note that the
natural projection $\widetilde{\text{pr}}_2:\mathcal
H_Y\to\oplus\bigoplus_{i=0}^{m-1}\widetilde{\mathcal H}_i$ is a
morphism of polarized Hodge structures.

It follows from the commutative diagram

\begin{picture}(280,64)
\put(156,48){\vector(-1,-1){18}} \put(158,50){$Y$}
\put(165,48){\vector(1,-1){18}}\put(142,42){$g$}
\put(174,42){$\sigma$} \put(130,20){$X$} \put(160,15){$f $}
\put(145,25){\vector(1,0){35}} \put(187,20){$V$}
\end{picture}
\newline
that $g^*: \mathcal H_X\to \mathcal H_Y$ is an embedding of
polarized Hodge structures, such that $g^*(f^*(\mathcal
H_V))=\sigma^*(\mathcal H_V)$ and
$g^*(M_V^\bot)\subset(\sigma^*(M_V))^\bot$. Therefore we can
identify $\mathcal H_X$ with its image $g^*(\mathcal H_X)$,
$f^*(\mathcal H_V)$ with  $\sigma^*(\mathcal H_V)$, and $M_V^\bot$
with $g^*(M_V^\bot)$. Under this identification, the projection
$\text{pr}^\bot$ is identified with the restriction of
$\widetilde{\text{pr}}_2$ to $g^*(\mathcal H_X)$. Therefore
$\text{pr}^\bot$ is a morphism of Hodge structures as the
composition of two morphisms of Hodge structures: namely, of the
embedding $g^*$ and the projection $\widetilde{\text{pr}}_2$. \qed \\

Let $\dim_{\mathbb C} X=2k$. Believing in Hodge Conjecture,  the
elements of $A_X=M_X\cap H^{k,k}$ will be called {\it algebraic}
and the module $T_X=\{ \gamma\in M_{X}\, \mid \, Q(\gamma,
\alpha)=0\, \, \text{for all}\, \, \alpha\in A_X\}$ will be called
the module of {\it transcendental} $n$-cycles on $X$, $n=2k$. It
is easy to see that the polarized Hodge structure $\mathcal H_X$
induces a polarized Hodge structure $\mathcal T_X=\{ T_X,
H_{T}^{p,q},Q_{T}\}$, where $H_T^{p,q}=(T\otimes \mathbb C)\cap
H_X^{p,q}$ and $Q_T$ is the restriction of $Q$ to $T_X$.

If $S$ is a smooth projective surface, then the form $Q=Q_S$ is
symmetric unimodular and, by Index Theorem, its signature is equal
to $(2h^{2,0}+1,h^{1,1}-1)$.   Note that the polarized Hodge
structure $\mathcal T_S$ on the transcendental cycles on a smooth
projective surface $S$ is a bi-rational invariant of the surface $S$. \\ \\
{\bf Non-decomposability Conjecture.} {\it  The polarized Hodge
structure $\mathcal T_S$ on the transcendental cycles on  a smooth
projective surface $S$  is non-decomposable.}
\\

Let $V\subset \mathbb P^5$ be a smooth cubic fourfold. It is known
(see \cite{T1} and \cite{T2}) that the moduli space of cubic
fourfolds contains several families of rational cubics.
Nevertheless, the following conjecture is well known.
\\ \\
{\bf Non-rationality Conjecture.} {\it A generic cubic fourfold is
non-ratio\-nal.} \\

The aim of the present note is to show that Non-rationality
Conjecture follows from  Non-decomposability Conjecture. To
formulate the precise statement, we fix one of smooth cubic
fourfolds $V_0$. For each cubic fourfold $V$ we can identify the
lattice $(M_V,Q_V)$ with the lattice $(M_{0}, Q)$, where
$M_{0}=H^4(V_0,\mathbb Z)$ and $Q=Q_{V_0}$. Let $\lambda =L^2\in
M_{0}$, where $L\in H^2(V_0,\mathbb Z)$ is the class of the
hyperplane section of $V_0\subset \mathbb P^5$. It is well known
(see, for example, \cite{K-K}) that each smooth cubic fourfold $V$
has the following Hodge numbers:   $h^{4,0}=h^{0,4}=0$,
$h^{3,1}=h^{1,3}=1$, and $h^{2,2}=21$. Consider the polarized
Hodge structures $\mathcal H_V=\{ M_V, H^{p,q},Q_V\}$ on the
fourth cohomology groups of the smooth cubic fourfolds $V$. Since
$h^{4,0}=0$ and $h^{3,1}=1$,  the polarized Hodge structure
$\mathcal H_V=\{ M_0, H^{p,q}_V,Q\}$ is defined by a non-zero
element $\omega \in H^{3,1}_V\subset M_0\otimes\mathbb C$. It
follows from Hodge -- Riemann bilinear relations that
$Q(\lambda,\omega)= Q(\omega,\omega)=0$ and
$Q(\omega,\overline{\omega})<0$. Therefore the classifying space
$D$ of polarized Hodge structures of the smooth cubic fourfolds
coincides with
$$D=\{ \omega \in \mathbb P^{22}\mid \,
Q(\omega,\lambda)=Q(\omega,{\omega})=0,\, \,
Q(\omega,\overline{\omega})<0\}.$$ The point  $\omega(V)\in D$,
corresponding to the polarized Hodge structure $\mathcal H_V=\{
M_0, H^{p,q}(V),Q\}$ of a smooth cubic fourfold $V$, is called the
{\it periods} of $V$. By \cite{G} and by Global Torelli Theorem
proved in \cite{V}, the set $D_0$ of periods of the smooth cubic
fourfolds is an open subset of $D$.

Let $P=\{ \mu\in M_0\, \mid Q(\lambda, \mu)=0\}$ be the submodule of
$M_0$ consisting  of the primitive elements. For each endomorphism
$A\in \text{End}(M_0)$ such that $A(P)\subset P$ and the restriction
$A_{\mid P}$ of $A$ to $P$ is not proportional to the identity
automorphism of $P$, denote by $$E_A=\{ \omega \in D\subset \mathbb
P^{22}\, \mid \omega \, \, \text{is an eigenvector of}\, \, A\} .$$
Obviously, $E_A$ is the intersection of $D$ and a finite number of
linear subspaces of $\mathbb P^{22}$ of codimension at least $2$.
Therefore $E_A$ is a proper closed analytic subvariety of $D$.  Put
$E=\cup E_A$, then $E$ is the union of countably many proper closed
analytic subvarieties of $D$. Therefore $D\setminus E$ is everywhere
dense in $D$.

For each $\mu \in M_0$, $\mu$ is not proportional to $\lambda$, we
put
$$B_{\mu}=\{ \omega \in D\subset \mathbb P^{22}\, \mid \, Q(\mu, \omega)=0\}$$ and
$B=\cup B_{\mu}$. Similarly to above, $B$ is the union of countably
many proper closed analytic subvarieties of $D$. Therefore
$D_0\setminus (E\cup B)$ is also everywhere dense in $D_0$.

\begin{prop} If Non-decomposability Conjecture is true,
then a smooth cubic fourfold $V$ is non-rational if its periods
$\omega(V)\in D_0\setminus (E\cup B)$.
\end{prop}

\proof Assume that a smooth cubic fourfold $V$ is rational and its
periods $\omega(V)\in D_0\setminus (E\cup B)$.

Note that the sublattice of transcendental elements $T_V$ of the
lattice $M_V$ coincides with the sublattice of primitive elements
$P$ if the periods $\omega(V)\in D_0\setminus B$.

Since $V$ is rational, then there is a bi-rational map $r:\mathbb
P^4 \dashrightarrow V$. By Hironaka Theorem, there is a composition
$\tau =\tau_n\circ\dots\circ\tau_1 :X\to \mathbb P^4$ of monoidal
transformations $\tau_i: X_i\to X_{i-1}$ with non-singular centers
$C_{i-1}\subset X_{i-1}$, where $X_0=\mathbb P^4$ and $X_n=X$, such
that $f=r\circ \tau :X\to V$ is a bi-rational morphism.

Each $\tau_i$ induces an inclusion of unimodular polarized Hodge
structures $\tau_i^*: \mathcal H_{X_{i-1}}\hookrightarrow \mathcal
H_{X_{i}}$ such that $\mathcal H_{X_{i}}\simeq \tau_i^*( \mathcal
H_{X_{i-1}})\oplus (\tau_i^*( \mathcal H_{X_{i-1}})^{\bot}$. In
particular,
$$H^4(X_{i},\mathbb Z))=(\tau_i^*( H^4(X_{i-1},\mathbb Z))\oplus
(\tau_i^*( H^4(X_{i-1},\mathbb Z))^{\bot}.$$ Moreover, if $\dim
C_{i-1}\leq 1$, then $(\tau_i^*( H^4(X_{i-1},\mathbb
Z))^{\bot}\subset H^{2,2}_{X_i}$, and if $C_{i-1}$ is a smooth
surface, then $(\tau_i^*( \mathcal H_{X_{i-1}})^{\bot}\simeq
-\mathcal H_{C_{i-1}}(-1)$, where $$-\mathcal
H_{C_{i-1}}(-1)=(M_{C_{i-1}},H_{i-1}^{p,q},-Q_{C_{i-1}})$$ and
$H_{i-1}^{p,q}=H_{C_{i-1}}^{p-1,q-1}$. Therefore $H^{3,1}_{X_i}=
\tau_i^*(H^{3,1}_{X_{i-1}})\oplus H^{3,1}_{i-1}$, $H^{4,0}_{X_i}=0$,
and $Q_{X_i}(u,v)=0$ for $u\in \tau_i^*(H^{3,1}_{X_{i-1}}\oplus
H^{1,3}_{X_{i-1}})$, $v\in H^{3,1}_{{i-1}}\oplus H^{1,3}_{{i-1}}$.

As a consequence, we obtain a decomposition of the unimodular
polarized Hodge structure
$$\mathcal H_X=\tau^*(\mathcal H_{\mathbb P^4})\oplus\bigoplus_{i=0}^{n-1}\mathcal
H_i$$ in the direct sum of polarized Hodge structures, where
$\mathcal H_i$ is the contribution in $\mathcal H_X$ of the
$(i+1)$th monoidal transformation. By induction, we have $\mathcal
H_i\simeq (\tau_i^*( \mathcal H_{X_{i-1}})^{\bot}$ and,
consequently,
$$\mathcal T_X\simeq -\oplus \mathcal T_{C_i}(-1),$$ where the sum is
taken over all surfaces $C_i$ with $p_g\geq 1$, since $H^4(\mathbb
P^4,\mathbb C)=H^{2,2}_{\mathbb P^4}$, $h^{2,2}(\mathbb P^4)=1$, and
$H^{2}(C_i,\mathbb C)=H^{1,1}_{C_i}$ if the geometric genus $p_g$ of
the surface $C_i$ is equal to zero.

The morphism $f$ induces a morphism of Hodge structures $f^*:
\mathcal H_V\to \mathcal H_X$. By Lemma 1,  the Hodge structure
$\mathcal H_X$ is decomposed in the direct sum $f^*(\mathcal
H_V)\oplus \mathcal H_V^\bot$ of the polarized Hodge structures.

\begin{lem} Let $V$ be a smooth cubic fourfold from Proposition 1
and let $\tau =\tau_n\circ\dots\circ\tau_1 :X\to \mathbb P^4$ and
$f=r\circ \tau :X\to V$ be morphisms described above. Then there is
$i_0$ such that the polarized Hodge structure $\mathcal T_{C_{i_0}}$
on the transcendental cycles of the surface $C_{i_0}$ is decomposed
in the direct sum of polarized Hodge structures $\mathcal
T^{\prime}_{C_{i_0}}\simeq -f^*(\mathcal T_V)(1)$ and 
$\mathcal T_{C_{i_0}}^{\prime\prime}$.
\end{lem}

\proof We have two decompositions
$$\mathcal H_X=f^*(\mathcal
H_V)\oplus \mathcal H_V^\bot=\tau^*(\mathcal H_{\mathbb
P^4})\oplus\bigoplus_{i=0}^{n-1}\mathcal H_i$$ of the polarized
Hodge structure $\mathcal H_X$. Consider a non-zero element
$\omega\in f^*(H^{3,1}_V)$. It can be represented in the form
$$\omega=\sum_{i=0}^{n-1}\omega_i,$$
where $\omega_i\in H^{3,1}_i$. Let $i_0$ be an index such that
$\omega_{i_0}\neq 0$. Put $\omega_{i_0}^\bot=\sum_{i\neq
i_0}\omega_i$ and let $\text{pr}_{i_0}:\mathcal H_X\to \mathcal
H_{i_0}$, $\text{pr}^{\bot}_{i_0}:\mathcal H_X\to \oplus_{i\neq i_0}
\mathcal H_{i}$,  $\text{pr}:\mathcal H_X\to f^*(\mathcal H_{V})$,
 $\text{pr}^{\bot}:\mathcal H_X\to \mathcal H_{V}^\bot$
be the natural projections. Note that all these projections are
defined over $\mathbb Z$. We have
$\text{pr}(\text{pr}_{i_0}(\omega))=a\omega$ for some $a\in\mathbb
C$, since $\dim f^*(H^{3,1}_V)=1$ and $\text{pr}$, $\text{pr}_{i_0}$
are morphisms of Hodge structures.

Let us show that $a=1$. Indeed, the restriction of
$\text{pr}\circ\text{pr}_{i_0}$ to $f^*(M_V)\simeq M_V$ induces an
endomorphism of $M_V$ such that
$$\text{pr}\circ\text{pr}_{i_0}(f^*(P))\subset f^*(P),$$ since
$\text{pr}\circ\text{pr}_{i_0}$ is a morphism of Hodge structures,
$M_V\cap H^{2,2}_V=\mathbb Z \lambda$ and $P=T_V$ by assumption.
Therefore $(\text{pr}\circ\text{pr}_{i_0})_{\mid
f^*(P)}=a\cdot\text{Id}$, since
$\text{pr}(\text{pr}_{i_0}(\omega))=a\omega$ (that is, $\omega$ is
an eigenvector of $(\text{pr}\circ\text{pr}_{i_0})_{\mid f^*(M_V)}$
and $a$ is its eigenvalue) and, by assumption, $\omega(V)\in
D_0\setminus (E\cup B)$. Therefore $a\in \mathbb Z$, since
$(\text{pr}\circ\text{pr}_{i_0})_{\mid f^*(P)}$ is an endomorphism
of $\mathbb Z$-module $f^*(P)\simeq P$.  Let
$\omega^\bot=\text{pr}^\bot(\omega_{i_0})\in H^{3,1}_X$. Then we
have
$$\begin{array}{ll}
\omega_{i_0}=& a\omega +\omega^\bot, \\
\omega^\bot_{i_0}= & (1-a)\omega -\omega^\bot . \end{array}$$ By
Hodge -- Riemann bilinear relations, $Q_X(\gamma,\overline
\gamma)\leq 0$ for $\gamma \in H^{3,1}_X$ and
$Q_X(\gamma,\overline\gamma)=0$ if and only if $\gamma =0$.
Therefore, without loss of generality, we can assume that
$Q_X(\omega,\overline \omega)=-1$ and $Q_X(\omega^\bot,\overline
\omega^\bot)=b\leq 0$.  We have $Q_X(\omega_{i_0},\overline
\omega_{i_0})=-a^2+b$ and, consequently, $a$ and $b$ can not
simultaneously  be equal to zero, since $\omega_{i_0}\neq 0$.
Besides, we have
$$Q_X(\omega,\omega^\bot)=Q_X(\omega,\overline
\omega^\bot)=Q_X(\overline \omega, \omega^\bot)=Q_X(\overline
\omega,\overline \omega^\bot)=Q_X(\omega_{i_0},\overline
\omega_{i_0}^\bot)=0.$$ Therefore
$$Q_X(\omega_{i_0},\overline
\omega_{i_0}^\bot)=Q_X(a\omega +\omega^\bot,(1-a)\overline\omega
-\overline\omega^\bot)=a(1-a)(-1)-b=0,$$ that is, $a^2-a-b=0$, and
hence
$$a=\frac{1\pm\sqrt{1+4b}}{2}.$$
But, $a\in \mathbb Z$ and $b\leq0$, therefore $b=0$, $a=1$ and hence
$\omega^\bot =0$, since $Q_X(\omega^\bot,\overline\omega^\bot)=b=0$.
Therefore, $\omega\in H^{3,1}_{i_0}$. Besides, we showed that
$(\text{pr}\circ\text{pr}_{i_0})_{\mid f^*(P)}=\text{Id}$.

Let us show that $f^*(T_V)\subset T_{i_0}$. As above, we have two
decompositions $$\mathcal T_X=f^*(\mathcal T_V)\oplus f^*(\mathcal
T_V)^\bot=\bigoplus_{i=0}^{n-1}\mathcal T_i$$ of the polarized Hodge
structure $\mathcal T_X$. Each $\gamma \in f^*(T_V)$ can be written
in the form $\gamma= \gamma_{i_0}+\gamma^\bot_{i_0}$, where
$\gamma_{i_0}=\text{pr}_{i_0}(\gamma)\in T_{i_0}$ and
$\gamma^\bot_{i_0}=\text{pr}_{i_0}^\bot(\gamma)\in \oplus_{i\neq
i_0} T_i$. Then, since $(\text{pr}\circ\text{pr}_{i_0})_{\mid
f^*(P)}=\text{Id}$ and $P=T_V$, we have
$$\begin{array}{ll}
\gamma_{i_0}= & \gamma +\gamma^\bot, \\  \gamma^\bot_{i_0}= &
-\gamma^\bot, \end{array}$$ where
$\gamma^\bot=\text{pr}^\bot(\gamma_{i_0})\in f^*(T_V)^\bot$. Denote
by $p_{(3,1)}: M_X\otimes C\to H^{3,1}_X$ the natural projection.
Then, by definition of transcendental cycles, $p_{(3,1)}(\gamma)\neq
0$ for each non-zero element $\gamma\in T_X$, since the Hodge number
$h^{4,0}(X)=0$. In particular, $p_{(3,1)}(\gamma)= a_{\gamma}\omega$
for some $a_{\gamma}\in \mathbb C$, $a_{\gamma}\neq 0$, since
$f^*(H_V^{3,1})=\mathbb C\omega$. Therefore
$p_{(3,1)}(\gamma_{i_0})= a_{\gamma}\omega$, since $\omega\in
f^*(H_V^{3,1})\subset H^{3,1}_{i_0}$, and hence
$p_{(3,1)}(\gamma^\bot_{i_0})=0$. From this, we have
$\gamma^\bot_{i_0}=0$, that is, $f^*(T_V)\subset T_{i_0}$.

As a consequence, we obtain that
\begin{equation} \label{decomp}
\mathcal T_{i_0}=f^*(\mathcal T_V)\oplus \mathcal
T_{i_0}^{\prime\prime}\end{equation} is the direct sum of
polarized Hodge structures, where $T_{i_0}^{\prime\prime}$ is a
polarized Hodge substructure of $f^*(\mathcal T_V)^\bot$. To
complete the proof of Lemma 2, recall that $\mathcal T_{i_0}\simeq
-\mathcal T_{C_{i_0}}(-1)$, where $\mathcal T_{C_{i_0}}$ is the
polarized Hodge structure on the transcendental $2$-cycles on the
smooth projective surface $C_{i_0}$. \qed
\begin{lem}  There does not exist a smooth
projective surface $S$ such that $\mathcal T_V\simeq -\mathcal
T_S(-1)$, where $V$ is a smooth cubic fourfold whose periods
$\omega(V)\in D_0\setminus B$.
\end{lem}
\proof If  $\mathcal T_V\simeq -\mathcal T_S(-1)$ for some surface
$S$, then its geometric genus $p_g=h^{2,0}(S)=h^{3,1}(V)=1$ and
$\text{rk}\, T_S=\text{rk}\, T_V=22$. Since $p_g=1$, the surface
$S$ is not a ruled surface. Note that the lattice $T_S$ is a
bi-rational invariant. Therefore we can assume that $S$ is a
minimal model and its second Betti number
$b_2(S)=h^{1,1}(S)+2h^{2,0}(S)\geq 23$, since $\text{rk}\,
T_S=\text{rk}\, T_V=22$ and $S$ should have also algebraic
$2$-cycles. Therefore $h^{1,1}(S)$ should be greater than $20$,
since $h^{2,0}(S)=1$.

On the other hand, it follows from the classification of algebraic
surfaces that $K^2_S\geq 0$, where $K_S$ is the canonical class of
$S$. Denote by $\chi(S)=1-h^{1,0}(S)+h^{2,0}(S)$ the algebraic Euler
characteristic of the surface $S$ and by
$e(S)=2+h^{1,1}_S+2h^{2,0}(S)-4h^{1,0}(S)$ its topological Euler
characteristic. Note that  $h^{1,0}(S)=\dim H^0(S,\Omega_S^1)\geq
0$. By Noether formula, we have $\chi(S)=\frac{K_S^2+e(S)}{12}$ and
hence
$$ h^{1,1}(S)=10+10h^{2,0}(S) -8h^{1,0}(S)-K^2_S=20 -8h^{1,0}(S)-K^2_S,$$
since $h^{2,0}(S)=1$. Therefore $ h^{1,1}(S)\leq 20$. \qed \\

It follows from Lemma 3 that in decomposition (\ref{decomp}), the
summand $\mathcal T_{i_0}^{\prime\prime}$ is non-trivial which
contradicts Non-decomposability Conjecture.

\ifx\undefined\bysame
\newcommand{\bysame}{\leavevmode\hbox to3em{\hrulefill}\,}
\fi

\end{document}